\documentclass[12pt]{article}
\usepackage{amssymb}
\usepackage{latexsym,bm}
\usepackage{graphicx}
\usepackage{amsmath}
\usepackage{mathrsfs}

\date{}
\begin{document}
\title{The minimum number of detours in graphs\footnote{E-mail address:
\tt zhan@math.ecnu.edu.cn}}
\author{\hskip -10mm Xingzhi Zhan\\
{\hskip -10mm \small Department of Mathematics, East China Normal University, Shanghai 200241, China}}\maketitle
\begin{abstract}
A longest path in a graph is called a detour. It is easy to see that a connected graph of minimum degree at least $2$ and order at least $4$
has at least $4$ detours. We prove that if the number of detours in such a graph of order at least $9$ is odd, then it is at least $9,$ and this lower bound
can be attained for every order. Thus the possibilities $3,$ $5$ and $7$ are excluded. Two open problems are posed.
\end{abstract}


\vskip 8mm

We consider finite simple graphs and use terminology and notations in [4] . Following Kapoor, Kronk, and Lick [3], we call a longest path in a graph $G$ a {\it detour} of $G.$ This concise
term has now been widely used (e.g. [1] and [2]). The {\it order} of a graph is its number of vertices. We denote by $V(G)$ and $E(G)$ the vertex set
 and edge set of a graph $G,$ respectively. For vertices $x$ and $y,$ an {\it $(x,y)$-path} is a path with endpoints $x$ and $y.$ We denote by $\delta(G)$ the minimum degree of a graph $G,$ and by $N(x)$
 the neighborhood of a vertex $x.$ If $u,\,v$ are two vertices on a path $P,$ then $P[u,v]$ denotes the subpath of $P$ with endpoints $u$ and $v.$ A basic fact about detours is that a detour of a
 connected graph $G$ of order $n$ has order at least ${\min}\{2\delta(G)+1,n\}.$

It is easy to see that a connected graph of minimum degree at least $2$ and order at least $4$ has at least $4$ detours. We will show that if the number of detours in such a graph
of order at least $9$ is odd, then it is at least $9,$ and this lower bound can be attained for every order. Thus the possibilities $3,$ $5$ and $7$ are excluded. At the end we pose two open problems.

{\bf Notation.} $f(G)$ denotes the number of detours in a graph $G.$

{\bf Theorem 1.}  {\it The minimum number of detours in a connected graph of minimum degree at least $2$ and order at least $4$ is $4.$}

{\bf Proof.} Let $G$ be a connected graph of order at least $4$ with $\delta(G)\ge 2$ and let $P:\, x_1,x_2,\ldots,x_k$ be a detour of $G.$
Then $N(x_1)\subseteq V(P)$ and $N(x_k)\subseteq V(P).$ Since $\delta(G)\ge 2,$ $x_1$ has a neighbor $x_i$ with $i\ge 3$ and $x_k$ has a neighbor $x_j$ with $j\le k-2.$
If $i=k$ or $j=1,$ then $G$ contains a $k$-cycle and we clearly have $f(G)\ge 4.$ Next suppose $3\le i\le k-1$ and $2\le j\le k-2.$

Case 1. $i\le j.$

$G$ has at least the following four detours:
$$
P,\quad P[x_1,x_j]\cup x_jx_k\cup P[x_k,x_{j+1}],
$$
$$
P[x_{i-1},x_1]\cup x_1x_i\cup P[x_i,x_k],\,\,\, P[x_{i-1},x_1]\cup x_1x_i\cup P[x_i,x_j]\cup x_jx_k\cup P[x_k,x_{j+1}].
$$
See Figure 1.
\vskip 3mm
\par
 \centerline{\includegraphics[width=3.3 in]{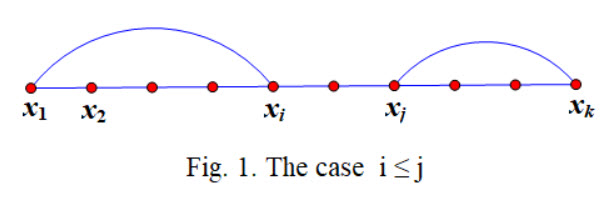}}
\par

Case 2. $i>j.$

$G$ has at least the following six detours:
$$
P,\,\,\, P[x_1,x_j]\cup x_jx_k\cup P[x_k,x_{j+1}],\,\,\, P[x_{i-1},x_1]\cup x_1x_i\cup P[x_i,x_k],
$$
$$
P[x_{j-1},x_1]\cup x_1x_i\cup P[x_i,x_k]\cup x_kx_j\cup P[x_j,x_{i-1}],\,\,\, P[x_{i+1},x_k]\cup x_kx_j\cup P[x_j,x_1]\cup x_1x_i\cup P[x_i,x_{j+1}],
$$
$$
P[x_{j-1},x_1]\cup x_1x_i\cup P[x_i,x_j]\cup x_jx_k\cup P[x_k,x_{i+1}].
$$
See Figure 2.
\vskip 3mm
\par
 \centerline{\includegraphics[width=3.3 in]{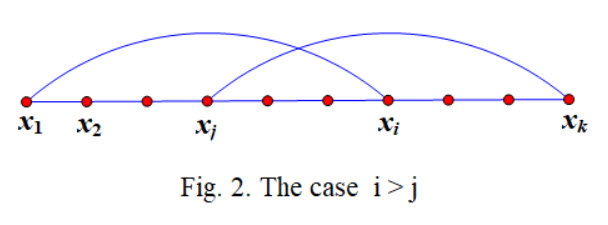}}
\par

This shows $f(G)\ge 4.$ Conversely, for every order $n\ge 4$ we construct a graph $G_n$ of order $n$ with $\delta(G)\ge 2$ satisfying
$f(G_n)=4.$ $G_4=C_4,$ the $4$-cycle. $G_5$ is the bowtie, the graph consisting of two triangles sharing one vertex. $G_6$ consists of a triangle and a $4$-cycle sharing one vertex.
$G_7,$ $G_8$ and $G_9$ are depicted in (a), (b) and (c) of Figure 3, respectively.
\vskip 3mm
\par
 \centerline{\includegraphics[width=4.4 in]{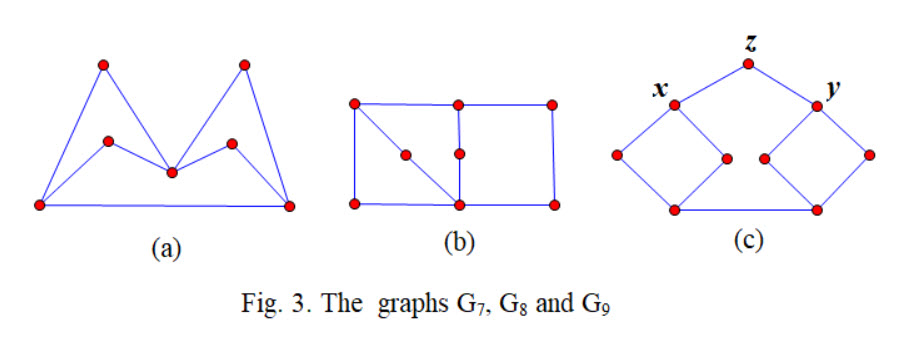}}
\par
For $n\ge 10,$ $G_n$ is obtained from $G_9$ in (c) of Figure 3 by replacing the path $x,z,y$ by an $(x,y)$-path of order $n-6.$ $\Box$

{\bf Remark 1.} Note that for $n\ge 7,$ the graphs $G_n$ in the above proof of Theorem 1 are 2-connected. Thus, if we replace ``minimum degree at least $2$" by ``2-connected"
in Theorem 1, we obtain the same conclusion for graphs of order at least $7.$

{\bf Remark 2.} In the above proof of Theorem 1, an edge $e$ of $P$ appears on at least four detours unless (1) $i\le j$ and  $e=x_{i-1}x_i$ or $e=x_jx_{j+1}$ or (2)
$i=j+1$ and $e=x_ix_j.$

{\bf Theorem 2.} {\it Let $G$ be a connected graph of minimum degree at least $2$ and order at least $9.$ If $f(G)$ is an odd number, then $f(G)\ge 9.$ Furthermore,
the lower bound $9$ can be attained for every order by both graphs of connectivity $1$ and graphs of connectivity $2.$}

{\bf Proof.} We first prove that  if $f(G)$ is an odd number, then $f(G)\ge 9.$ It suffices to show that either $f(G)\ge 8$ or $f(G)=4$ or $f(G)=6.$

We make the following conventions: (1) For a positive integer $r,$ ``$r$ detours" means ``$r$ pair-wise distinct detours"; (2) for an edge $e$ of $G$ and a detour $D,$ we say that
$e$ appears on $D$ if $e\in E(D).$

Let $P:\, x_1,x_2,\ldots,x_k$ be a detour of $G.$ If there is another detour $Q$ with $V(Q)\not=V(P),$ by the proof of Theorem 1, there are at least $4$ detours with the same vertex set $V(P)$
and there are at least $4$ detours with the same vertex set $V(Q).$ These detours are clearly distinct. Hence we have $f(G)\ge 8.$ Next suppose that all detours of $G$ have $V(P)$ as
their vertex set.

Recall that an edge $e$ of $G$ is called a {\it chord} of a path $R$ if the two endpoints of $e$ lie in $R$ but $e\not\in E(R).$ A chord $e$ of $R$ is called an {\it inner chord} if both
endpoints of $e$ are internal vertices of $R.$ Otherwise $e$ is called a {\it boundary chord.} A detour $D$ is called a {\it basic detour} if no inner chord of $P$ is an edge of $D;$ otherwise
$D$ is called a {\it non-basic detour}.

Let the order of $G$ be $n.$ If $G$ is hamiltonian, then $f(G)\ge n\ge 9.$ Next assume that $G$ is non-hamiltonian.

Since $P$ is a detour, $N(x_1)\subseteq V(P)$ and $N(x_k)\subseteq V(P).$ The condition $\delta(G)\ge 2$ implies that $x_1$ has a neighbor $x_i$ with $i\ge 3$ and $x_k$ has a neighbor $x_j$ with $j\le k-2.$
If $i=k$ or $j=1,$ then $G$ has a $k$-cycle $C$ which contains $P.$ Since $P$ is a detour, $C$ must be a Hamilton cycle, contradicting our assumption that $G$ is non-hamiltonian.
Hence $3\le i\le k-1$ and $2\le j\le k-2.$ We distinguish two cases.

{\bf Case 1.} Every detour of $G$ is a basic detour.

We need consider only the boundary chords of $P.$

{\bf Subcase 1.1.} $P$ contains exactly two boundary chords.

As analyzed in the proof of Theorem 1, in this case $f(G)=4$ or $f(G)=6.$

{\bf Subcase 1.2.} $P$ contains exactly three boundary chords.

Without loss of generality, let $x_1x_q$ be the third chord of $P$ with $q\not=i.$ Note that the two boundary chords $x_1x_i$ and $x_1x_q$ are in symmetric positions.
If $q>i$ we may interchange the roles of $x_1x_i$ and $x_1x_q.$ Thus we may and do assume that $q<i.$

Suppose $i\le j.$ We have four basic detours not containing the edge $x_1x_q.$ If $3\le q\le i-2,$ we have exactly the following two detours containing the edge $x_1x_q:$
$$
P[x_{q-1},x_1]\cup x_1x_q\cup P[x_q,x_k],\,\,\,P[x_{q-1},x_1]\cup x_1x_q\cup P[x_q,x_j]\cup x_jx_k\cup P[x_k,x_{j+1}].
$$
Hence $f(G)=6.$ If $q=i-1,$ we have exactly the following four detours containing the edge $x_1x_q:$
$$
P[x_2,x_q]\cup x_qx_1\cup x_1x_i\cup P[x_i,x_k],\,\,\,P[x_2,x_q]\cup x_qx_1\cup x_1x_i\cup P[x_i,x_j]\cup x_jx_k\cup P[x_k,x_{j+1}],
$$
$$
P[x_{q-1},x_1]\cup x_1x_q\cup P[x_q,x_k],\,\,\,P[x_{q-1},x_1]\cup x_1x_q\cup P[x_q,x_j]\cup x_jx_k\cup P[x_k,x_{j+1}].
$$
Hence $f(G)=8.$

Suppose $i>j.$ We have six basic detours not containing the edge $x_1x_q.$ In this case it is easy to check that there are at least two detours containing the edge $x_1x_q$
by considering the subgraph $P\cup x_1x_q\cup x_kx_j.$ Thus $f(G)\ge 8.$

{\bf Subcase 1.3.} $P$ contains at least four boundary chords.

Based on Subcase 1.2, we deduce that $f(G)\ge 8$ in this case.

{\bf Case 2.} $G$ contains a non-basic detour.

Claim 1. Every edge in a detour appears on at least two detours.

This claim can be verified by checking the proof of Theorem 1, replacing $P$ there by the detour in question. In fact,
except for possible one or two edges, every edge in a detour appears on at least four distinct detours. See Remark 2 above.

Since $G$ contains a non-basic detour, some inner chord $e$ of $P$ is an edge of a detour. By Claim 1, there are at least two detours containing $e$ as an edge.
Thus $G$ has at least two non-basic detours. If $i>j,$ we have six basic detours, and consequently $f(G)\ge 8.$
By Subcases 1.2 and 1.3, if $P$ contains at least three boundary chords, then $G$ contains at least six basic detours. Again we have $f(G)\ge 8.$
If an inner chord of $P$ appears on at least four detours, then we have at least four non-basic detours. It follows that $f(G)\ge 8.$

It remains to treat the case when (1) $i\le j,$ (2) $P$ contains exactly two boundary chords and (3) every inner chord of $P$ appears on
at most three detours. Next we make these three assumptions.

Let $D:\, y_1,y_2,\ldots,y_k$ be a detour of $G.$ Suppose $y_c$ is a neighbor of $y_1$ and $y_d$ is a neighbor of $y_k.$ As in the proof of Theorem 1,
there are four detours (if $c\le d$) or six detours (if $c>d$) whose edges belong to $E(D)\cup \{y_1y_c,\,y_ky_d\}.$ We denote by $\Psi(D)$ the set of these four
or six detours according as $c\le d$ or $c>d.$ When we write $\Psi (D)$ we assume that the two boundary chords $y_1y_c$ and $y_ky_d$ have been prescribed.

Claim 2. If an inner chord $h$ of $P$ appears on a detour $D$ such that $\Psi(D)\cap\Psi(P)\not=\phi,$ then one of the two endpoints of $h$
belongs to the set $\{x_{i-1}, x_{j+1}\}.$

Let $D=y_1,y_2,\ldots,y_k$ with $h\in E(D).$ Suppose $y_c$ is a neighbor of $y_1$ and $y_d$ is a neighbor of $y_k.$
Since $h$ appears on at most three detours, by Remark 2 after the proof of Theorem 1, if $c>d$ we must have $c=d+1$ and then $D$ is contained in a cycle
which must be a Hamilton cycle since $D$ is a detour, contradicting our assumption that $G$ is non-hamiltonian. Thus $c\le d$ and then either $h=y_{c-1}y_c$ or $h=y_dy_{d+1}.$
Note that since $h$ is an inner chord of $P$, the two endpoints of $h$ cannot be $x_1$ or $x_k.$ Let $R\in\Psi(D)\cap\Psi(P).$ Then $R$ does not contain $h.$
Each of the two detours in $\Psi(D)$ not containing $h$ has one endpoint which is an endpoint of $h.$ Thus one endpoint $v$ of $R$ is an endpoint of $h.$
Since the endpoints of the four detours in $\Psi(P)$ are $x_1,x_k,x_{i-1},x_{j+1},$ we deduce that $v\in \{ x_1,x_k,x_{i-1},x_{j+1}\}$ but $v\not\in \{ x_1,x_k\}.$
Hence $v\in \{x_{i-1},x_{j+1}\}$.

{\bf Subcase 2.1.} $G$ has a detour which contains at least two inner chords of $P.$

Let $h$ and $e$ be two inner chords of $P$ that appear on one common detour. Consider the subgraph $G^{\prime}=P\cup x_1x_i\cup x_kx_j\cup h\cup e.$
The path $T=P[x_{i-1},x_1]\cup x_1x_i\cup P[x_i,x_j]\cup x_jx_k\cup P[x_k,x_{j+1}]$ is a detour of $G^{\prime}$ with endpoints $x_{i-1}$ and $x_{j+1},$ which is also
a detour of $G.$ Since we have assumed that $G$ is non-hamiltonian,  $x_{i-1}$ and $x_{j+1}$ are non-adjacent. By Claim 2, each of $h$ and $e$ has exactly one endpoint
in the set $\{x_{i-1},x_{j+1}\}.$ Now in $G^{\prime},$ the detour $T$ has four boundary chords $x_{i-1}x_i,$ $x_{j+1}x_j,$ $h$ and $e.$
By Subcase 1.3 above (replacing $P$ there by $T$), we obtain $f(G)\ge f(G^{\prime})\ge 8.$

{\bf Subcase 2.2.} Every non-basic detour contains exactly one inner chord of $P.$

Denote by $\Omega$ the set of the inner chords of $P$ that appear on at least one detour.
By the above Claim 1, if one inner chord of $P$ appears on a detour, then there are at least two detours containing that chord. Thus, if $|\Omega|\ge 2$
then we have at least four non-basic detours, and consequently we have $f(G)\ge 8.$ Next suppose $|\Omega|=1$ and let $\Omega=\{x_sx_t\}$ with $2\le s \le t-2.$
Recall that we have assumed $i\le j.$ Using Claim 2, we deduce that $f(G)=6$ if (1) $t=i-1;$ (2) $s=j+1;$ (3) $s=i-1$ and $i+2\le t\le j;$ (4) $i\le s\le j-2$ and $t=j+1.$
In all other cases $f(G)\ge 8.$ This completes the proof that  if $f(G)$ is an odd number, then $f(G)\ge 9.$

Next for every integer $n\ge 9$ we construct a graph $H_n$ of order $n$ and connectivity $1$ which contains exactly $9$ detours. Every $H_n$ is traceable.
$H_9,H_{10}$ and $H_{11}$ are depicted in Figure 4.
\vskip 3mm
\par
 \centerline{\includegraphics[width=5.7 in]{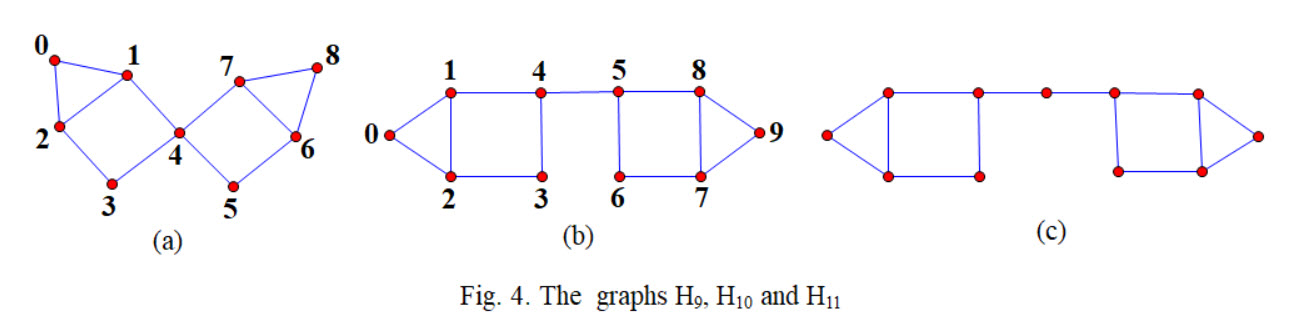}}
\par
For $n\ge 11,$ $H_n$ is obtained from $H_{10}$ by subdividing the edge $(4,5)$ $n-10$ times, i.e., replacing the edge $(4,5)$ by a $(4,5)$-path of order $n-8.$
The $9$ detours in $H_9$ are
$$
(0,1,2,3,4,5,6,7,8), \,\,\, (0,1,2,3,4,5,6,8,7),\,\,\, (0,1,2,3,4,7,8,6,5),
$$
$$
(1,0,2,3,4,5,6,7,8), \,\,\, (1,0,2,3,4,5,6,8,7),\,\,\, (1,0,2,3,4,7,8,6,5),
$$
$$
(3,2,0,1,4,5,6,7,8), \,\,\, (3,2,0,1,4,5,6,8,7),\,\,\, (3,2,0,1,4,7,8,6,5).
$$

Finally for every integer $n\ge 9$ we construct a graph $M_n$ of order $n$ and connectivity $2$ which contains exactly $9$ detours. Every $M_n$ is traceable.
$M_9$ and $M_{10}$ are depicted in Figure 5.
\vskip 3mm
\par
 \centerline{\includegraphics[width=3.3 in]{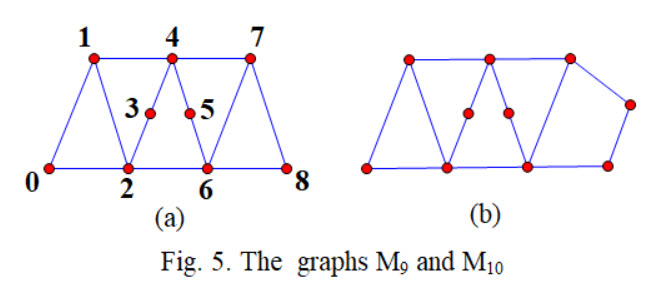}}
\par
For $n\ge 10,$ $M_n$ is obtained from $M_9$ by subdividing the edge $(7,8)$ $n-9$ times.
Observe that  $M_9$ is obtained from $H_9$ in Figure 4 (a) by adding the edge $(2,6),$ and any detour of $M_9$ cannot contain the edge $(2,6).$
Hence $M_9$ and $H_9$ have the same set of detours, in particular, the same number of detours, i.e., $9.$
Note that each detour of $M_9$ contains the edge $(7,8).$ Thus for every $n\ge 10,$ $M_n$ has the same number of detours as $M_9.$ $\Box$

Finally we pose two problems. Recall that $f(G)$ denotes the number of detours in a graph $G.$

{\bf Problem 1.} Let $k$ and $n$ be integers with $3\le k\le n-2.$ Denote by $\Gamma (k,n)$ the set of connected graphs with minimum degree $k$ and order $n.$ 
Define
$$
a(k,n)={\rm min}\{f(G)|\, G\in \Gamma (k,n)\}.
$$
Determine $a(k,n).$

{\bf Problem 2.} Let $k,$ $n$ and $\Gamma(k,n)$ be as in Problem 1. Define
$$
b(k,n)={\rm min}\{f(G)|\, G\in \Gamma (k,n)\,\,\,{\rm and}\,\,\,f(G)\,\,\,{\rm is}\,\,\,{\rm odd}\}.
$$
Determine $b(k,n).$

Perhaps for sufficiently large orders $n,$ $a(k,n)$ and $b(k,n)$ are independent of $n.$

\vskip 5mm
{\bf Acknowledgement.} This research  was supported by the NSFC grant 12271170 and Science and Technology Commission of Shanghai Municipality
 grant 22DZ2229014.


\begin{thebibliography}{99}
\bibitem{1} L.W. Beineke, J.E. Dunbar and M. Frick, Detour-saturated graphs, J. Graph Theory 49(2005), 116–134.
\bibitem{2} G. Chartrand, G.L. Johns and S.L. Tian, Detour distance in graphs, Quo vadis, graph theory?, 127–136, Ann. Discrete Math., 55, North-Holland, Amsterdam, 1993.
\bibitem{3} S.F. Kapoor, H.V. Kronk and D.R. Lick, On detours in graphs, Canad. Math. Bull. 11(1968), 195–201.
\bibitem{4} D.B. West, Introduction to Graph Theory, Prentice Hall, Inc., 1996.
\end{thebibliography}
\end{document}